\documentclass[12pt]{article}
\usepackage{amsmath,amssymb,color,amscd}

\def\seq#1#2#3{#1_{#2},\,\ldots,#1_{#3}}

\def\vv{{\underline{v}}}

\def\nunu{\underline{\nu}}
\def\tt{{\underline{t}}}

\newcommand{\Var}{\ensuremath{\mathcal{V}_{\mathbb{C}}}}
\def\1{\underline{1}}
\def\R{\mathbb R}
\def\P{\mathbb P}

\def\L{\mathbb L}
\def\LL{{\cal L}}

\def\A{\mathbb A}
\def\F{\mathbb F}
\def\Z{\mathbb Z}

\def\C{\mathbb C}

\def\OO{{\cal O}}

\newenvironment{definition}
{\smallskip\noindent{\bf Definition\/}:}{\smallskip\par}

\newenvironment{remark}
{\smallskip\noindent{\bf Remark\/}.}{\smallskip\par}

\title{On Poincar\'e series of filtrations
\footnote{Math. Subject Class. 16W70, 13A18, 16W60.
Keywords: filtrations, Poincar\'e series, valuations.}
}

\author{
A.~Campillo,
\and F.~Delgado,\thanks{Supported by the grants
MTM2012-36917-C03-01 / 02
(both grants with the help of FEDER Program).} \and S.M.~Gusein-Zade
\thanks{
Supported by the
grants RFBR--13-01-00755, NSh--5138.2014.1.
} }

\date{}
\begin{document}
\def\eps{\varepsilon}

\maketitle

\begin{abstract}
In this survey one discusses the notion of the Poincar\'e series of multi-index filtrations,
an alternative approach to the definition, a method
of computation of the Poincar\'e series based on the notion of integration with respect to the Euler
characteristic (or rather on an infinite-dimensional version of it), generalizations of the notion
of the multi-variable Poincar\'e series based on the notion of the motivic integration, and relations
of the latter ones with some zeta functions over finite fields and with generating series of
Heegaard-Floer homologies of algebraic links.
\end{abstract}

\section*{Introduction}\label{sec0}

Let $(V,0)$ be a germ of a complex analytic space. A decreasing {\em filtration} on the ring $\OO_{V,0}$
of germs of functions on $(V,0)$ is a decreasing sequence of vector subspaces (sometimes ideals) $\{J_i\}$
in $\OO_{V,0}$:
$$
\OO_{V,0}=J_0\supset J_1 \supset J_2\supset\ldots
$$
An important invariant of a filtration is its {\em Poincar\'e series}
\begin{equation}\label{eq1}
P(t)=\sum_{i=0}^{\infty} \dim(J_i/J_{i+1})\cdot t^i\,.
\end{equation}
The Poincar\'e series $P(t)$ is defined if all the factors $J_i/J_{i+1}$ are finite-dimensional or,
in other words, if $\OO_{V,0}/J_i$ is finite-dimensional for each $i$.
It is a formal power series with integer coefficients, that is an element of the ring $\Z[[t]]$ of formal
power series in the variable $t$.

As an example one can consider the filtration defined on the ring $\OO_{\C^n,0}$ of germs of functions in
$n$ variables by an irreducible germ of a curve $(C,0)\subset(\C^n,0)$. Let $\varphi:(\C,0)\to(\C^n,0)$ be
a parameterization (an uniformization) of the curve $(C,0)$, i.e. a complex analytic map such that
${\rm Im}\,\varphi=C$ and $\varphi$ is an isomorphism between $\C$ and $C$ outside the origin.
For a function germ $g:(\C^n,0)\to \C$, let $\nu(g)$ be the order of zero of the composition
$g\circ\varphi:(\C,0)\to\C$ at the origin, i.e. $g\circ\varphi(\tau)=a(g)\tau^{\nu(g)}+{\rm term\ of\ higher\ order}$,
where $a(g)\ne 0$. If $g\circ\varphi\equiv 0$, one assumes that $\nu(g)=+\infty$. The function
$\nu:\OO_{\C^n,0}\to\Z\cup\{+\infty\}$ is a {\em valuation} on $\OO_{\C^n,0}$,
i.e. $\nu(\lambda g)=\nu(g)$ for a non-zero complex number $\lambda$, $\nu(g_1+g_2)\ge\min(\nu(g_1),\nu(g_2))$,
and $\nu(g_1g_2)=\nu(g_1)+\nu(g_2)$.
The filtration corresponding to the valuation $\nu$ is defined by
$$
J_i=\{g\in \OO_{\C^n,0}: \nu(g)\ge i\}\,.
$$
(The subspace $J_i$ is an ideal in $\OO_{\C^n,0}$.)

Assume that $n=2$, i.e. $(C,0)$ is an irreducible plane curve singularity, and let $(C,0)$ be defined by
an equation $f=0$, where $f:(\C^2,0)\to(\C,0)$ is a germ of a function in two variables.
In \cite{FAA} it was shown that the Poincar\'e series of the filtration
$\{J_i\}$ (being a priori an analytic invariant of $(C,0)$) coincides with the so called monodromy zeta function
(see the definition, e.g., in \cite{AGV}) of the germ $f$. This coincidence was obtained by direct
computations of the both objects in the same terms. Up to now there is no conceptual proof of this relation.

Generalization of the notion of the Poincar\'e series of a usual (one-index) filtration to multi-index filtrations
is far from being straightforward. (Even there are different versions of a definition of a multi-index filtration
which demand somewhat different definitions of the Poincar\'e series.)
There were found generalizations of the described relation between the Poincar\'e series and the monodromy
zeta function in the multi-index setting.

In this survey we discuss the notion of the Poincar\'e series of multi-index filtrations, a method
of computation of the Poincar\'e series based on the notion of integration with respect to the Euler
characteristic (or rather on an infinite-dimensional version of it), generalizations of the notion
of the multi-variable Poincar\'e series based on the notion of the motivic integration, and relations
of the latter ones with some zeta functions over finite fields and with generating series of
Heegaard-Floer homologies of algebraic links.

\section{Multi-index filtrations}\label{sec1}

A function $\nu:\OO_{V,0}\to \Z_{\ge0}\cup \{+\infty\}$ is a {\em valuation} if
\begin{enumerate}
 \item[1)] $\nu(\lambda g)=\nu(g)$ for $\lambda\in\C$, $\lambda\ne 0$;
 \item[2)] $\nu(g_1+g_2)\ge\min(\nu(g_1),\nu(g_2))$;
 \item[2)] $\nu(g_1g_2)=\nu(g_1)+\nu(g_2)$.
\end{enumerate}
A function $\nu:\OO_{V,0}\to \Z_{\ge0}\cup \{+\infty\}$ is called an {\em order function} if it possesses the
properties 1) and 2) above.

A (usual, one-index) filtration $\{J_i\}$ on $\OO_{V,0}$ can be defined by the order function
$$
\nu(g):=\max\{i:g\in J_i\}\,.
$$
Namely, $J_i=\{g\in \OO_{V,0}: \nu(g)\ge i\}$. The order function $\nu$ will be called {\em finite}
if the quotients $\OO_{V,0}/J_i$ are finite dimensional for all $i$.

There are two possible generalizations of the notion of a filtration to a multi-index setting, i.e.
when the elements of the filtration are numbered not by non-negative integers, but by elements
of $\Z_{\ge0}^r=\{\vv=(v_1, \ldots, v_r): v_i\in\Z_{\ge0}\}$. For $\vv'=(v'_1, \ldots, v'_r)$
and $\vv''=(v''_1, \ldots, v''_r)$, one says that $\vv'\ge\vv''$ if $v'_i\ge v''_i$ for all $i$.
A multi-index ($r$-index) filtration on $\OO_{V,0}$ is defined by a system $\{J_{\vv}\}$ of
vector subspaces of $\OO_{V,0}$ for all $\vv\in\Z_{\ge0}^r$. In any case one assumes that
$J_{\underline{0}}=\OO_{V,0}$ ($\underline{0}=(0,\ldots,0)$) and
$J_{\vv'}\subset J_{\vv''}$ for $\vv'\ge\vv''$.
In the first (less restrictive) approach these are the only conditions. In the second one we should assume that
the filtration is defined by a collection $\{\nu_i\}_{1=1, \ldots, r}$ of order functions on $\OO_{V,0}$.
This means that
\begin{equation}\label{eq2}
J_{\vv}=\{g\in\OO_{V,0}: \nunu(g)\ge\vv\}\,,
\end{equation}
where $\nunu(g):=(\nu_1(g), \ldots, \nu_r(g))$. In terms of the subspaces $J_{\vv}$ themselves this condition
can be formulated in the following way.
For two elements $\vv'=(v'_1, \ldots, v'_r)$ and $\vv''=(v''_1, \ldots, v''_r)$ from $\Z_{\ge0}^r$, let us define
$\max(\vv', \vv'')\in\Z_{\ge0}^r$ as $(\max(v'_1, v''_1), \ldots, \max(v'_r, v''_r))$. One can show that a
multi-index filtration $J_{\vv}$ is defined by a collection of order functions if and only if
$J_{\max(\vv', \vv'')}=J_{\vv'}\cap J_{\vv''}$ for any $\vv'$ and $\vv''$. In what follows
(except Section \ref{sec2})
we shall consider multi-index filtrations defined by collections of order functions.

An example of a filtration not defined by order functions can be constructed in the following way (see \cite{L}).
Let $I$ be an ideal in $\OO_{V,0}$ and let $(W,0)$ be a subscheme of $(V,0)$ defined by the ideal $I$.
This means that $\OO_{W,0}=\OO_{V,0}/I$. Let $\{J_{\vv}\}_{\vv\in\Z^r_{\ge0}}$ be an arbitrary multi-index
filtration on $\OO_{V,0}$
(possibly defined by a collection of order functions). The {\em embedded filtration} on $\OO_{W,0}$ induced by
the filtration $\{J_{\vv}\}$
is the filtration $\{I_{\vv}\}$ defined by
$$
I_{\vv}=\phi(J_{\vv})\,,
$$
where $\phi:\OO_{V,0}\to\OO_{W,0}=\OO_{V,0}/I$ is the canonical factorization map. One can show that in general
the filtration $\{I_{\vv}\}$ on $\OO_{W,0}$ is not defined by order functions.

\section{Poincar\'e series of multi-index filtrations}\label{sec2}

The notion of the Poincar\'e series of a multi-index filtration was first introduced in \cite{CDK}
(for filtrations defined by collections of valuations). The definition is far from being
a straightforward generalization of (\ref{eq1}).

Let $\{J_{\vv}\}$ be a multi-index filtration on $\OO_{V,0}$. We shall assume that the subspaces $J_{\vv}$
are defined not only for $\vv\in\Z^r_{\ge0}$, but for all $\vv\in\Z^r$. This can be done using the definition
$$
J_{\vv}:=J_{\max(\vv,\underline{0})}\quad{\rm for\ \ }\vv\in\Z^r\,,
$$
where $\underline{0}=(0,\ldots, 0)$. If the filtration $\{J_{\vv}\}$ is defined by a collection
$\{\nu_i\}_{i=1, \ldots, r}$ of order functions on $\OO_{V,0}$, this means that $J_{\vv}$ is given by
the equation (\ref{eq2}) for all $\vv\in\Z^r$.

Let
$$
\LL(t_1, \ldots, t_r):=
\sum\limits_{\vv\in\Z^r} \dim(J(\vv)/J(\vv+\1))\tt^{\,\vv}\,.
$$
where $\1=(1, 1,\ldots, 1)\in\Z^r$, $\tt:=(t_1, \ldots, t_r)$, $\tt^{\,\vv}:=t_1^{v_1}\cdot\ldots\cdot t_r^{v_r}$.
This equation makes sense if and only if all the quotients $J(\vv)/J(\vv+\1)$ are finite-dimensional.
This is equivalent to the condition that all the quotients $\OO_{V,0}/J(\vv)$ are finite-dimensional.
The series $\LL(\tt)$ is an element of the set (a lattice, that is a free $\Z$-module) of formal Laurent series in
$t_1$, \dots, $t_r$ with integer coefficients. This set is not a ring: one cannot, in general, multiply two series
of this sort. However it is a module over the ring $\Z[t_1,\ldots,t_r]$ of polynomials in $t_1$, \dots, $t_r$ with
integer coefficients. (In fact also a module over the ring of Laurent polynomials in $t_1$, \dots, $t_r$.)
If $v_i<0$ at least for one $i$, one has $J_{\vv}=J_{(v_1, \ldots, v_{i-1},0,v_{i+1},\ldots, v_r)}$.
This implies that $\LL(t_1, \ldots, t_r)\cdot\prod_{i=1}^r(t_i-1)$ does not contain monomials with negative exponents,
that is it is an element of the ring $\Z[[t_1, \ldots, t_r]]$ of formal power series in $t_1$, \dots, $t_r$.
In \cite{CDK} the authors offered the following notion.

\begin{definition}
 The {\em Poincar\'e series} of the filtration $\{J_{\vv}\}$ is
\begin{equation}\label{Poincare}
 P(t_1, \ldots, t_r)=\frac{\LL(t_1, \ldots, t_r)\cdot\prod_{i=1}^r(t_i-1)}{t_1\cdot\ldots\cdot t_r-1}.
\end{equation}
\end{definition}

This definition can be (and sometimes is) used for the both types of filtrations described in Section~\ref{sec1}.
However, properties of the
Poincar\'e series and methods of their computation are quite different in these cases. In order to indicate
the filtration or (in the case if the filtration is defined by order functions) the system of order functions,
one uses the notation $P_{\{J_{\vv}\}}(\tt)$ or $P_{\{\nu_i\}}(\tt)$.

By computing the right hand side of (\ref{Poincare}) and taking into account that $J_{\vv}=J_{\underline{0}}=\OO_{V,0}$ for $\vv\le\underline{0}$,
it is easy to check that for every $\vv$ the coefficient at $\tt^{\vv}$ in the Poincar\'e series $P(t_1, \ldots, t_r)$
is equal to
$$
-\sum_{K\subset\{1,\ldots, r\}}(-1)^{\vert K\vert}\dim \OO_{V,0}/J_{\vv+\underline{1}_K}=-\chi(C'_{\bullet,\vv})\,,
$$
where
$\underline{1}_K$ is the element of the lattice $\Z^r$ whose components with the numbers from $K$
are equal to $1$ and the other components are equal to $0$,
$C'_{\bullet,\vv}$ is the chain complex of vector spaces whose chain spaces are
$$
C'_{i,\vv}=\bigoplus_{K\subset\{1,\ldots,r\}, \vert K\vert=i}\ \OO_{V,0}/J_{\vv+\underline{1}_K}
$$
and the chain map $\partial_i:C'_{i,\vv}\to C'_{i-1,\vv}$ on the components corresponding to $K$ and $K'$ with $\vert K\vert=i$
and $\vert K'\vert=i-1$ is equal to 0 if $K'\not\subset K$ and is equal to $(-1)^{\ell}$ times the natural map
$\OO_{V,0}/J_{\vv+\underline{1}_K}\to \OO_{V,0}/J_{\vv+\underline{1}_{K'}}$ if $K'=K\setminus\{k\}$ and $k$ is the $\ell$th
integer in the natural order among those in $K$.

Since $J_{\vv+\underline{1}}\subset J_{\vv+\underline{1}_K}$ for each $K$, one can also define the chain complexes
$C_{\bullet,\vv}$ and $C''_{\bullet,\vv}$ given by
\begin{eqnarray*}
C_{i,\vv}&=&\bigoplus_{K\subset\{1,\ldots,r\}, \vert K\vert=i}\ J_{\vv+\underline{1}_K}/J_{\vv+\underline{1}}\\
\mbox{and}&{\ }&{\ }\\
C''_{i,\vv}&=&\bigoplus_{K\subset\{1,\ldots,r\}, \vert K\vert=i}\ \OO_{V,0}/J_{\vv+\underline{1}}
\end{eqnarray*}
respectively with the chain maps defined in the same way as above. Since the vector spaces in the sum
for $C''_{i,\vv}$ are all the same,
one has $\chi(C''_{\bullet,\vv})=0$ (since $C''_{\bullet,\vv}$ is nothing but the complex counting
the augmented homologies
of the $(r-1)$-dimensional simplex with the coefficients in $\OO_{V,0}/J_{\vv+\underline{1}}$).
Thus one gets that the coefficient
at $\tt^{\vv}$ in $P(t_1, \ldots, t_r)$ is also given by $\chi(C_{\bullet,\vv})=-\chi(C'_{\bullet,\vv})$ and therefore one has the equations
\begin{equation}\label{Antonio1}
P(t_1, \ldots, t_r)= \sum_{\vv\in\Z_{\ge 0}}\chi(C_{\bullet,\vv})\tt^{\vv}\,=\,-\sum_{\vv\in\Z_{\ge 0}}\chi(C'_{\bullet,\vv})\tt^{\vv}\,.
\end{equation}
Moreover, one also has the equations involving the homology spaces:
\begin{equation}\label{Antonio2}
\chi(C_{\bullet,\vv})=\sum_{i=0}^r (-1)^i h_{i,\vv}\,,\quad \chi(C'_{\bullet,\vv})=\sum_{i=0}^r (-1)^i h'_{i,\vv}\,,
\end{equation}
where $h_{i,\vv}=\dim H_i(C_{\bullet,\vv})$, $h'_{i,\vv}=\dim H_i(C'_{\bullet,\vv})$. One has
$h_{i,\vv}=h'_{i-1,\vv}$ for all $\vv$ and $i$.

The chain complexes considered above are defined even if the spaces $J_{\vv}/J_{\vv+\underline{1}}$ or, equivalently,
the spaces $\OO_{V,0}/J_{\vv}$ are not finite dimensional. In fact it is possible that the values $h_{i,\vv}$
(or equivalently the values $h'_{i,\vv}$) are finite for all $\vv$ and $i$ even if the spaces
$J_{\vv}/J_{\vv+\underline{1}}$ are not finite dimensional. This motivates the
following (more universal) definition of the Poincar\'e series of a filtration.

\begin{definition}
If $h_{i,\vv}$ is finite for all $\vv\ge 0$ and $i=0,1,\ldots, r$ then the {\em Poincar\'e series}
of the filtration $\{J_{\vv}\}$ is defined by
\begin{equation}\label{Antonio3}
P_{\{J_{\vv}\}}=\sum_{\vv\in\Z_{\ge 0}^r}\sum_{i=0}^r (-1)^i h_{i,\vv}\tt^{\vv}
=-\sum_{\vv\in\Z_{\ge 0}^r}\sum_{i=0}^r (-1)^i h'_{i,\vv}\tt^{\vv}\,.
\end{equation}
\end{definition}

For filtrations given by order functions one of which is finite, in \cite{CL} it is proved
that all the dimensions
$h_{i,\vv}$ are finite and therefore the Poincar\'e series can be defined by the equation (\ref{Antonio3}).
In particular, this is the case if one of them is the order function defined by the maximal ideal
$\mathfrak{m}$ of $\OO_{V,0}$ or if it is a valuation centred at this ideal.

When all the order functions are not finite, the Poincar\'e series is also defined in some cases.
For instance, this happens if $(V,0)$ is the germ of a toric variety, the order functions from the
collection are monomial valuations and the vector spaces
$$
H_{0,\vv}=J_{\vv}/\left(J_{\vv+\underline{1}_1}+J_{\vv+\underline{1}_2}+\ldots+J_{\vv+\underline{1}_r}\right)
$$
are finite dimensional: see \cite{CL}. In this case one can prove that
$\chi(C_{\bullet,\vv})=h_{0,\vv}=\dim H_{0,\vv}$.
As an example, if $(V,0)=(\C^2,0)$, $r=2$, and $\nu_1$ and $\nu_2$ are the $x$-order and the $y$-order valuations
respectively, one gets $P_{\{\nu_1,\nu_2\}}(t_1,t_2)=\frac{1}{(1-t_1)(1-t_2)}$.

Another example is when $(V,0)$ is irreducible and one considers a collection $\{\nu_1,\ldots,\nu_r\}$
of discrete valuations on it and a function germ $f\in\OO_{V,0}$ such that $\nu_i(f)\le +\infty$. One
has two filtrations: the one
on $\OO_{V,0}$ defined by the collection $\{\nu_1,\ldots,\nu_r\}$ and the embedded filtration
$\{I_{\vv}\}$ on the hypersurface
germ $\{f=0\}\subset(V,0)$. Assume that the Poincar\'e series $P_{\{\nu_i\}}(\tt)$ for the initial filtration
on $\OO_{V,0}$ can be defined by the equation (\ref{Antonio3}). In \cite{L} and \cite{CL} it is shown
that the Poincar\'e series
$P_{\{I_{\vv}\}}(\tt)$ of the embedded filtration can be defined by the same equation and one has
$$
P_{\{I_{\vv}\}}(\tt)=(1-\tt^{\underline{q}})P_{\{\nu_i\}}(\tt)\,,
$$
where $\underline{q}=(\nu_1(f), \ldots, \nu_r(f))$.

\section{Poincar\'e series and the Hilbert functions}\label{sec3}

A somewhat more traditional invariant of a multi-index filtration is the so-called Hilbert function
$h:\Z_{\ge 0}^r\to\Z_{\ge 0}$ defined by
\begin{equation}\label{Hilbert}
 h(\vv)=\dim \OO_{X,0}/J_{\vv}\,.
\end{equation}
This function is defined if and only if all the quotients $\OO_{X,0}/J_{\vv}$ have finite dimensions,
i.e. if and only if the Poincar\'e series can be defined by the equation (\ref{Poincare}).
Note that the definition (\ref{Hilbert}) makes sense for all $\vv\in\Z^r$. It can be described by its
generating function which can be defined either as a formal power series
$$
\widetilde{H}(\tt)=\sum_{\vv\in\Z^r_{\ge0}}h(\vv)\cdot\tt^{\vv}\,,
$$
or as a formal Laurent series
$$
H(\tt)=\sum_{\vv\in\Z^r}h(\vv)\cdot\tt^{\vv}\,.
$$
One can easily see that the series $\widetilde{H}(\tt)$ and $H(\tt)$ determine each other.

One can show that
$$
P(\tt)=-\tt^{-\1}H(\tt)\prod_{i=1}^r(t_i-1)\,.
$$
(This equation can be deduced, e.g., from (\ref{Antonio1}).)
This means that the Poincar\'e series of a multi-index filtration is determined by the generating function
$H(\tt)$, that is by the Hilbert function.

Since $t_i-1$ is a divisor of zero in the module of formal Laurent series
($(t_i-1)(\ldots+t_i^{-2}+t_i^{-1}+1+t_i^1+t_i^2+\ldots)=0$), one cannot, in general, restore
the Hilbert function from the Poincar\'e series and thus the Hilbert function is a more fine invariant.

As an example one can take two germs of curves
$C= C_1\cup C_2\subset (\C^5,0)$ (respectively $C' = C_1'\cup C'_2\subset (\C^6,0)$) whose
branches $C_1$ and $C_2$ (resp. $C'_1$ and $C'_2$) are defined by
the parameterizations ($t\in \C$, $u\in \C$):
\begin{eqnarray*}
C_1 &= \{(t^2,t^3,t^2,t^4,t^5)\}\,,\quad
C_2 &= \{(u^2,u^3, u^4, u^2, u^6)\}  \\
(\mbox{resp.}\;  C'_1 &=  \{(t^3,t^4,t^5, t^4,t^5, t^6)\}\,,\quad
C'_2 &= \{(u^3, u^4,  u^5, u^5, u^6, u^7)\} \; ).
\end{eqnarray*}
The local rings $\OO_{C,0}$ and $\OO_{C',0}$ of the curves $C$ and
$C'$, as subrings of
the normalizations $\C\{t\}\times \C\{u\}$, are:
\begin{eqnarray*}
\OO_{C,0} & = & \C\{\,(t^2,u^2), (t^3, u^3), (t^2,u^4), (t^4, u^2), (t^5,
u^6)\,\}, \\
\OO_{C',0} & = & \C\{\,(t^3,u^3), (t^4,u^4), (t^5, u^5), (t^4,u^5),
(t^5, u^6), (t^6, u^7)\,\}\; .
\end{eqnarray*}
An easy computation shows that in both cases the Poincar\'e
series is the polynomial $P(\tt) = 1 + t_1^{3}t_2^3$. However
one has $h(3,3)=1$ for the curve $C'$ and $h(3,3)=3$ for $C$.

In \cite{Monats2014} it was shown that the Hilbert function contains the same information about
a filtration as the so-called generalized Poincar\'e series discussed in Section~\ref{sec5}.
For filtrations on the ring $\OO_{\C^2,0}$ of germs of functions in two variables defined by valuations
it was shown (\cite{Monats2014}) that the Hilbert function determines the topological type of the set of valuations,
that is the topological type of its minimal resolution. (In \cite{FAOM} it was shown that in this situation
the Poincar\'e series, in general, does not determine the topological type of the set of valuations.)

\section{Poincar\'e series of filtrations and integration with respect to the Euler characteristic}\label{sec4}

The definition (\ref{Poincare}) (and even the definition (\ref{eq1}) for a one-index filtration) does not give,
in general, a direct way to compute the Poincar\'e series of a filtration.

One method to compute the Poincar\'e series is based on a reformulation of the definition (\ref{Poincare})
in terms of the Euler characteristic.
Let $\{J_{\vv}\}$, $\vv\in\Z^r_{\ge0}$, be a multi-index filtration on $\OO_{V,0}$ defined by
a collection $\{\nu_i\}$ of order functions.
For $\vv\in\Z^r_{\ge0}$ and $g\in\OO_{V,0}$, $\nunu(g)=\vv$ if and only if
$g\in J_{\vv}$ and $g\notin J_{\vv+\1_i}$ for all $i=1,\ldots,r$, where $\1_i=(0,\ldots, 0,1,0,\ldots, 0)$
with $1$ at the $i$th place. Let
$$
F_{\vv}:=(J_{\vv}/J_{\vv+\1})\setminus\left(\bigcup_{i=1}^r J_{\vv+\1_i}/J_{\vv+\1}\right)
\subset J_{\vv}/J_{\vv+\1}\,.
$$
The set $F_{\vv}$ is the complement to an arrangement of vector subspaces ($J_{\vv+\1_i}/J_{\vv+\1}$, $i=1,\ldots, r$)
in a vector space ($J_{\vv}/J_{\vv+\1}$).
The disjoint union $\widehat{S}=\coprod_{\vv\in\Z_{\ge 0}^r}F_{\vv}$ of the sets $F_{\vv}$
is called the {\em extended semigroup} of the filtration $\{J_{\vv}\}$ or of the collection
of order functions $\{\nu_i\}$. It is really a semigroup if $\nu_i$ are valuations.
The set $F_{\vv}$ is invariant with respect to multiplication by non-zero complex numbers.
Let $\P F_{\vv}=F_{\vv}/\C^*\subset\P (J_{\vv}/J_{\vv+\1})$ be the projectivization of the set $F_{\vv}$.
It is the complement to an arrangement of projective subspaces in a projective space.

For a (good enough) topological space $X$ (say, for a semi-algebraic set) let $\chi(X)$ be its Euler characteristic
defined by:
\begin{equation}\label{Euler}
 \chi(X)=\sum_{q=1}^{\infty} (-1)^q \dim H_c^q(X;\R)\,,
\end{equation}
where $H_c^q(\bullet)$ are cohomology groups with compact support. Euler characteristic defined this way is not
a homotopy invariant. It is an invariant of the homotopy type defined in terms of proper maps (i.e. maps such that
the preimage of a compact subspace is compact). An important property of this Euler characteristic
is its additivity:
$\chi(X)=\chi(Y)+\chi(X\setminus Y)$ for a closed (semi-algebraic) subset $Y\subset X$. The Euler characteristic
defined in the same way through the usual cohomology groups is a homotopy invariant, but is not additive. For example,
for $X$ and $Y$ being the circle and its point respectively, the both Euler characteristics of $X$ and $Y$ are equal
to $0$ and $1$ respectively. (These two Euler characteristics coincide for compact spaces.) The complement
$X\setminus Y$ is homeomorphic to the open interval. Its ``usual'' Euler characteristic
(that is the Euler characteristic defined through the usual cohomology groups) is equal to $1$
(since the interval is homotopy equivalent to the point),
whence its Euler characteristic defined by (\ref{Euler}) is equal to $(-1)$. The additivity of the
Euler characteristic (\ref{Euler}) implies the inclusion-exclusion formula for it.

\begin{remark}
In fact these two Euler characteristics also coincide for complex quasi-projective varieties.
(A {\em quasi-projective variety} is a projective variety minus another projective variety.)
\end{remark}

Let the filtration $\{J_{\vv}\}$ be defined by a system $\{\nu_i\}$ of order functions, $i=1,\ldots, r$.
One can show that in this case
\begin{equation}\label{extended}
 P(\tt)=\sum_{\vv\in\Z^r_{\ge0}}\chi(\P F_{\vv})\tt^{\vv}\,.
\end{equation}
This follows essentially from two facts. First, the multiplication by $\prod_{i=1}^r(t_i-1)$ can be interpreted
as the inclusion-exclusion formula applied to the coefficients of a series. Second, the Euler characteristic
of the projectivization of a (finite-dimensional) complex vector space is equal to its dimension. Another way to
prove (\ref{extended}) is to show directly that $\chi(\P F_{\vv})=\chi(C_{\bullet,\vv})$ (see (\ref{Antonio1})).

The right hand side of the equation (\ref{Poincare}) and (\ref{extended}) can be interpreted as certain integrals with respect
to the Euler characteristic. This notion was first formulated precisely in \cite{Viro}. Here we shall formulate it
for the category of real constructible sets. A {\em semi-algebraic} real {\em set} is a subset of
a real projective space defined in affine charts by a finite collection of algebraic equations and
algebraic inequalities ($>$, $<$, $\ge$, $\le$). A {\em constructible} set is the union of a finite number
of semi-algebraic sets. The constructible subsets of a projective space constitute an algebra of sets.
The Euler characteristic of a semi-algebraic set is defined by the equation (\ref{Euler}).
The additivity of the Euler characteristic permits to extend it to all constructible sets.

Let $X$ be a constructible set and let $A$ be an abelian group. A function $\psi:X\to A$ is called
{\em constructible} if it has a finite number of values and for each $a\in A$ the level set $\psi^{-1}(a)$
is constructible. The additivity of the Euler characteristic on the algebra of constructible subsets of $X$
permits to use it as a sort of a (non-positive) measure for an analogue of integration. The {\em integral}
of a constructible function $\psi$ over $X$ {\em with respect to the Euler characteristic} is defined by
\begin{equation}\label{integral}
 \int_X\psi d\chi=\sum_{a\in A}\chi(\psi^{-1})a\,.
\end{equation}
The right hand side of this equation can be regarded as an integral sum. It makes sense since it contains finitely many
summands.

The integral with respect to the Euler characteristic possesses many properties of the usual integral.
In particular, one has the Fubini formula for it. Let $p:X\to Y$ be a constructible map. (This means that
its graph in $X\times Y$ is constructible.) Let us consider the function
$$
\int_{p^{-1}(y)}\psi d\chi
$$
on $Y$. One can show that this function is constructible (first it was proved in \cite{Var}) and one has
$$
\int_X\psi d\chi = \int_Y\left(\int_{p^{-1}(y)}\psi d\chi\right) d\chi\,.
$$

In these terms the equation (\ref{extended}) can be written as an integral with respect to the Euler characteristic:
\begin{equation}\label{integral_semigroup}
P(\tt)=\int_{\P\widehat{S}}\tt^{\vv} d\chi\,,
\end{equation}
where $\vv:\P\widehat{S}\to\Z_{\ge0}^r$ is the tautological function on $\P\widehat{S}$ sending the component
$\P F_{\vv}$ to $\vv$.

In \cite{IJM} the notion of the integral with respect to the Euler characteristic was extended to integrals
over the infinite-dimensional space $\P\OO_{V,0}$. The idea was to consider ``approximations'' of $\P\OO_{V,0}$
by the finite-dimensional projective spaces $\P J^k_{V,0}$, where $J^k_{V,0}:=\OO_{V,0}/\mathfrak{m}^{k+1}$
is the space of $k$-jets of functions on $(V,0)$. (Here $\mathfrak{m}$ is the {\em maximal ideal} in $\OO_{V,0}$,
i.e., the set of functions with the value $0$ at the origin.) After that the integral over the infinite-dimensional
space $\P\OO_{V,0}$ is defined as the limit of appropriate integrals over the finite-dimensional
spaces $\P J^N_{k,0}$.

The precise definition is the following. Let $d_k=\dim J^k_{V,0}=\mbox{codim\,}\mathfrak{m}^{k+1}$ and
let $\P^* J^k_{V,0}$ be the disjoint union of $\P J^k_{V,0}$
with a point (in some sense $\P^* J^k_{V,0}=J^k_{V,0}/\C^*$). One has natural maps
$\pi_k: \P{\OO}_{V,0} \to \P^* J^k_{V,0}$ and
$\pi_{k,\ell}: \P^* J^k_{V,0} \to \P^* J^{\ell}_{V,0}$ for $k \ge \ell$.
Over $\P J^\ell_{V,0} \subset \P^* J^\ell_{V,0}$ the map $\pi_{k,\ell}$
is a locally trivial fibration, the fibre of which
is a complex vector space of dimension $d_k-d_{\ell}$.
(E.g., for $V=\C^n$ this dimension is equal to ${n+k\choose k}-{n+\ell\choose \ell}$.)
A subset $X\subset \P{\cal O}_{\C^n,0}$ is called {\em cylindric} if
$X=\pi_k^{-1}(Y)$ for a constructible subset
$Y\subset \P J^k_{V,0} \subset \P^* J^k_{V,0}$.
For a cylindric subset $X\subset \P{\cal O}_{V,0}$ ($X=\pi_k^{-1}(Y)$,
$Y\subset \P J^k_{V,0}$) its {\em Euler characteristic} $\chi(X)$ is defined
as the Euler characteristic $\chi(Y)$ of the set $Y$.
A function $\psi: \P{\OO}_{V,0} \to A$ with values in
an abelian group $A$ with not more than countably many values
is called {\em cylindric} if, for each $a\ne 0$,
the set $\psi^{-1}(a)\subset \P{\OO}_{V,0}$ is cylindric.
Now the integral of a cylindric function $\psi$ over the space
$\P{\OO}_{V,0}$ with respect to the Euler characteristic is
defined by the same equation (\ref{integral}) (where $X=\P{\cal O}_{V,0}$)
with the only difference that the right hand side of it may contain infinitely many
summands and thus may make no sense in $A$ (``not to be convergent'').
If the right hand side of (\ref{integral}) makes sense in $A$, the function $\psi$ is said to be integrable.

This definition (together with the interpretation of the dimensions of vector spaces as the Euler
characteristics of their projectivizations) permits to rewrite the equation (\ref{extended}) as
\begin{equation}\label{Poincare-integral}
 P(\tt)=\int_{\P\OO_{V,0}}\tt^{\nunu(g)} d\chi
\end{equation}
(where one assumes that $t_i^{+\infty}=0$).

The equation (\ref{Poincare-integral}) turned out to give a powerful method to compute the Poincar\'e series
of filtrations in a number of cases: see, e.g., \cite{IJM}, \cite{DG}, \cite{Invent}. An example of computations
can be found in \cite[Theorem 6]{RMS}.

\section{Generalized Poincar\'e series and some their applications}\label{sec5}

The additivity of the Euler characteristic is the main property which permits to use it as a sort of a
measure. There are some other invariants of quasi-projective varieties possessing this property. As an
example one can take the Hodge-Deligne polynomial of a variety.
The most general additive invariant is the class of the variety in the corresponding Grothendieck ring (see below).
Therefore it can be considered as a universal (generalized) Euler characteristic and can be used for
integration. This permits to define more fine invariants of filtrations similar to the Poincar\'e series.

Let $K_0(\Var)$ be the Grothendieck ring of
quasi-projective varieties. It is generated by the classes $[X]$ of
such varieties subject to the relations:\newline
1) if $X_1\cong X_2$, then $[X_1]=[X_2]$;\newline
2) if $Y$ is Zariski closed in $X$, then $[X]=[Y]+[X\setminus Y]$.\newline
The multiplication in $K^0(\Var)$ is defined by the Cartesian product. The class $[X]\in K_0(\Var)$ of a
variety $X$ can be regarded as the generalized Euler characteristic $\chi_g(X)$ of it.
Let $\L$ be the class $[{\mathbb A}^1_{\C}]$ of the complex affine line.
The class $\L$ is not equal to zero in the ring $K_0(\Var)$.
Moreover the natural ring homomorphism
$\Z[x]\to K_0(\Var)$ which sends $x$ to $\L$ is an
inclusion. (This follows, e.g., form the fact that the Hodge-Deligne polynomial of $P(\L)$, $P$ is a
polynomial, is $P(uv)$.)
Let
$K_0(\Var)_{(\L)}$
be the localization of the
Grothendieck ring $K_0(\Var)$ by the class $\L$.
The natural homomorphism
$\Z[x]_{(x)}\to K_0(\Var)_{(\L)}$ is an inclusion as
well. A remarkable fact is that the class of the complex line $\L$ is a zero divisor in
$K_0(\Var)$ (see \cite{Borisov}). This fact does not affect the following
constructions and definitions.

\medskip

Let the filtration $\{J_{\vv}\}$ on the ring $\OO_{V,0}$ be defined by a system $\{\nu_i\}$ of order
functions, $i=1,\ldots, r$. The equations (\ref{extended}) and (\ref{integral_semigroup}) for the Poincar\'e
series
suggest the following definition
for the {\it generalized semigroup Poincar\'e series\/}
\begin{equation}\label{gen_extended}
{\widehat P}_g(\tt)
=\sum_{\vv\in\Z^r_{\ge0}}[\P F_{\vv}]\cdot\tt^{\,\vv}
= \int_{\P\hat{S}}\tt^{\nunu(g)} d\chi_g
\in K_0(\Var)[[t_1,\ldots, t_r]]\,.
\end{equation}

All connected components of $\P\widehat S$ (i.e.
projectivizations $\P F_\vv$ of the fibres $F_\vv$\,) are
complements to arrangements of projective subspaces in {\it
finite dimensional} (!) projective spaces. Because of that all the coefficients of the series
${\widehat P}_g(\tt)$ are polynomials in $\L$. Therefore we can write
${\widehat P}_g(t_1,\ldots,t_r)$ as a series ${\widehat P}_g(t_1,\ldots,t_r,\L)\in
\Z[[t_1,\ldots,t_r,\L]]$ in $t_1$, \dots, $t_r$, and $\L$.
One has $P(t_1,\ldots,t_r)={\widehat P}_g(t_1,\ldots,t_r,1)$.

In a similar way as the one used with the (non generalized) Poincar\'e series in the
previous section, one gets the following equation:
$$
\widehat P_g(\tt)
=\frac{\left(\sum\limits_{\vv\in\Z^r}[\P(J_{\vv}/J_{\vv+1})]\cdot\tt^{\,\vv}\right)\cdot
\prod\limits_{k=1}^r(t_k-1)}{t_1\cdot\ldots\cdot t_r -1}\ .
$$

\medspace

One can extend the notion of integration with respect to the generalized Euler characteristic to the
infinite dimensional space $\P\OO_{V,0}$ in a way similar to the one for the usual Euler characteristic.

\begin{definition}
For a cylindric subset $X\subset \P\OO_{V,0}$
($X=\pi_k^{-1}(Y)$,
$Y\subset \P J^k_{V,0}$\,, $Y$ is constructible),
its {\em generalized Euler characteristic} $\chi_g(X)$
is the element $[Y]\cdot\L^{-d_k}$ of the  ring
$K_0(\Var)_{(\L)}$, where $d_k = \dim (\OO_{V,0}/\mathfrak{m}_{V.0}^{k+1})$.
\end{definition}

The generalized Euler characteristic $\chi_g(X)$ is well defined
since, if $X = \pi^{-1}_{\ell}(Y')$,
$Y'\subset \P J^{\ell}_{V,0}$\,, $k\ge \ell$, then $Y$ is  a
Zariski locally trivial fibration over $Y'$ with the fibre
$\C^{d_k-d_\ell}$ and therefore $[Y]=[Y']\cdot
\L^{d_k-d_\ell}$.

\begin{definition}
Let $\psi: \P\OO_{V,0} \to A$ be a function with values in
an abelian group $A$ with countably many values.
The integral of $\psi$ over the space
$\P\OO_{V,0}$ with respect to
the generalized Euler characteristic
is
$$
\int_{\P\OO_{V,0}}\psi d\chi_{g} :=
\sum_{a\in A, a\ne 0} \chi_{g}(\psi^{-1}(a))\cdot a,
$$
where $\chi_{g}$ is the
generalized Euler characteristic
if this sum makes sense in
$A\otimes_{\Z}K_0(\Var)_{(\L)}$. If the integral
exists
(makes sense) the function $\psi$ is said to be integrable.
\end{definition}

\begin{definition}
The {\it generalized Poincar\'e series \/} of the collection $\{\nu_i\}$ of order functions
is
$$
P_g(\tt):= \int_{\P\OO_{V,0}} \tt^{\nunu} d\chi_g \, .
$$
\end{definition}

The subset of the projectivization $\P\OO_{V,0}$ where
$\tt^{\,\nunu(g)}$ is equal to $\tt^{\,\vv}$ (i.e., where $\nunu(g)=\vv$)
is the projectivization of the space
$J_{\vv}\setminus\bigcup\limits_{k=1}^{r} J_{\vv + \1_{\{k\}}}$.
Because of that all the coefficients of the series
$P_g(\tt)$ are polynomials in $\L^{-1}$. Therefore we shall write
$P_g(t_1,\ldots,t_r)$ as a series $P_g(t_1,\ldots,t_r,q)\in
\Z[[t_1,\ldots,t_r,q]]$ in $t_1$, \dots, $t_r$, and $q=\L^{-1}$.
One has $P(t_1,\ldots,t_r)=P_g(t_1,\ldots,t_r,1)$.

The dimension of the quotient of two vector subspaces $E' \subset E$ of $\OO_{V,0}$ can be computed
as $\chi(\P E\setminus \P E')$.
Substituting the usual Euler characteristic $\chi$ by the
generalized one $\chi_g$, one gets the following ``motivic version"
of the dimension: $``\dim_g"(E/E')=\chi_g(\P E \setminus\P E')$.
Let $\mbox{codim }E=a$, $\mbox{codim }E'=a'$ ($a' > a$). Then
$\chi_g(\P E \setminus\P E')=q^{a+1}+q^{a+2}+\ldots+q^{a'}=q^{a+1}\cdot\frac{1-q^{\,a'-a}}{1-q}$.
Applying this computation to
$J_{\vv}$ and $J_{\vv + \1}$, one gets
$$
\chi_g(\P J_{\vv}\setminus \P J_{\vv + \1} ) =
q^{\,h(\vv)+1}\cdot\frac{1-q^{\,h(\vv+\1)-h(\vv)}}{1-q}\cdot\tt^{\,\vv}\,.
$$
Therefore, if one defines
$$
L_g(\tt,q):=\sum_{\vv\in\Z^r}
q^{\,h(\vv)+1}\cdot\frac{1-q^{\,h(\vv+\1)-h(\vv)}}{1-q}\cdot\tt^{\,\vv}\,
$$
one can show that
$$
P_g(\tt,q)=\frac{L_g(\tt,q)\cdot\prod\limits_{k=1}^r
(t_k-1)}{t_1\cdot\ldots\cdot t_r-1}\,.
$$

Notice that the generalized Poincar\'e series and the generalized semigroup Poincar\'e
series are different. In fact they belong to different rings, but both specialize to
the Poincar\'e series $P(\tt)$.

\bigskip

Let $k$ be a perfect field.
The definition of the generalized Poincar\'e series can be extended
to the category of varieties defined over the field $k$. In this case one must consider
the Grothendieck ring $K_0(\mathcal{V}_{k})$ of reduced quasi-projective varieties
defined over the field $k$ and its localization
$K_0(\mathcal{V}_{k})_{(\L)}$ on the class of the affine line,
$\L = [\A^1_{k}]$.
This construction, in the particular case of curves defined over a finite field
$\F_{p^s}=k$, allows to compute some interesting invariants of
the local ring be means of the integration with respect to the generalized Euler
characteristic.

\medskip
Let $C$ be a complete geometrically irreducible, algebraic curve defined over the
field $\F_{p^s}$ and  let $\OO$ be the local ring of $C$ in a closed point $P\in C$.
Let $\pi: \widetilde{C}\to C$ be the normalization,  let $\seq{P}1r$ be the points on
$\widetilde{C}$ lying over $P$ and let $\OO_i$ ($1\le i\le r$) be the corresponding local rings of
$\widetilde{C}$ at the points $P_i$.
The integral closure of $\OO$ in its function field $K$ is just
the intersection $\OO_1\cap \ldots \cap \OO_r$, each $\OO_i$ is the valuation ring
of a valuation $\nu_i$ over $K$.
Let us denote
$\rho:=[\OO/\mathfrak{m}: \F_{p^s}]$ and
for $i=1,\ldots, m$,  $d_i:= [\OO_i/\mathfrak{m}_i: \F_{p^s}]$.

The valuations $\seq{\nu}1r$ define, as usual, the filtration
$$
J_{\vv} := \{g\in \OO | \nu_i(g)\ge v_i, 1\le i \le r\} \mbox{ for }
\vv\in\Z^r
$$
and so one has the corresponding generalized Poincar\'e series:
\begin{equation}
P_g(\tt, \L) = \int_{\P\OO} \tt^{\nunu(g)}d \chi_g \in
K_0(\mathcal{V}_{\F_{p^s}})_{(\L)}[[\seq t1r]]\; .
\end{equation}

St\"ohr in \cite{Stohr} introduces
the zeta function now called the St\"ohr one:
$$
\zeta (\OO, z):= \sum\limits_{\mathfrak{a} \supseteq \OO} \sharp \left
(\mathfrak{a}/ \OO \right )^{-z},
$$
where the summation runs through all fractional ideals containing
$\OO$,  $z \in \C$ with $\mathrm{Re}(z)>0$.
Putting $t=q^{-z}$ and writing $Z(\OO,t)$ instead of
$\zeta (\OO,z)$,
the St\"ohr zeta function splits into a
finite sum of partial zeta functions
$Z(\OO,t)=\sum_{(\mathfrak{b})}Z(\OO,\mathfrak{b},t)$, where
$\mathfrak{b}$ is a fractional $\OO$-ideal satisfying
$\overline{\OO} \cdot \mathfrak{b} = \overline{\OO} $. The
summation runs over a complete system of representatives of the
ideal class semigroup of $\OO$. For each partial zeta function one
has:
\begin{equation} \label{eq:Z_o_b_t}
Z(\OO,\mathfrak{b},t)=\sum_{\mathfrak{a} \supseteq \OO \; , \; \mathfrak{a} \sim
\mathfrak{b}} t^{\dim_{k} \left ( \mathfrak{a} / \OO \right )}.
\end{equation}
The notation $\mathfrak{a} \sim \mathfrak{b}$ means that $\mathfrak{a} = a^{-1}
\mathfrak{b}$ for a non-zero divisor $a \in K$.

In \cite{D_Moyano} it is shown that the partial zeta functions $Z(\OO,\mathfrak{b},t)$ can be computed
in terms of integrals with respect to the Euler characteristic:
$$
Z(\OO, \mathfrak{b},t)= \frac{(p^s-1)(p^s)^{\rho} t^{\deg
(\mathfrak{b})}}{((p^s)^{\rho}-1) (U_{\mathfrak{b}}:U_{\OO}) }
\int_{\mathbb{P}\mathfrak{b}} (p^s t)^{\underline{\nu}(g) \cdot
\underline{d}} d \chi.
$$
Here $\nunu(g)\cdot\underline{d}:= d_1\nu_1(g)+\cdots+d_r\nu_r(g)$,
$(U_{\mathfrak{b}}: U_{\OO})$ is the index of the subgroup $U_{\OO}$ of units of $\OO$ over the units of
the fractional ideal $\mathfrak{b}$,
and the degree $\deg(\mathfrak{b})$ of $\mathfrak{b}$ is the function such that
$\deg(\OO)=0$ and $\deg(\mathfrak{a}/\mathfrak{b})=\deg(\mathfrak{a})-\deg(\mathfrak{b})$.
Finally, $\chi(X) = \# \chi_g(X)$ is just the cardinality of the generalized Euler characteristic.

For $\mathfrak{b}=\OO$ one has
\begin{equation*}
Z(\OO,\OO,t) = \sum_{i=0}^{\infty} \# \{\mbox{principal ideals of } \OO \mbox{ of
codimension }i\}\cdot t^i \; .
\end{equation*}
In \cite{D_Moyano}, it is proved that
\begin{equation}\label{zeta1}
Z(\OO,\OO,t) = \frac{(p^s)^\rho (p^s-1)}{(p^s)^\rho-1} P_g((p^s t)^{d_1},\ldots,
(p^s t)^{d_r}; p^{s})\; .
\end{equation}
Notice that the right hand side of the equation [\ref{zeta1}]
has a simple expression in the totally rational case,
i.e. when $d_i=1$ for $i=1,\ldots, r$ (this implies also that $\rho=1$):
\begin{equation*}
Z(\OO,\OO,t) =  p^s P_g(p^s t,\ldots,
p^s t; p^{s})\; .
\end{equation*}

Also in this context,
Moyano and Z\'u\~niga-Galindo
(\cite{M_ZG})
showed the rationality of the generalized Poincar\'e series
$P_g(\tt, \L)$ and its relation with a zeta function
associated to the effective Cartier divisors on the curve $C$.

\medskip

One can say that the generalized Poincar\'e series is (in the sense of M.Kho\-va\-nov) a ``categorification''
of the usual Poincar\'e series. This means that the generalized Poincar\'e series is a series in one more
variable $q=\L^{-1}$ and its Euler characteristic (obtained by the substitution $q=1$) coincides with the
the usual Poincar\'e series. There are several ``categorifications'' of the Alexander polynomial of a
knot or a link. Thus one may hope that the generalized Poincar\'e series can coincide with (or be related to)
one of these ``categorifications''. One of them called the Heegaard-Floer link homology was constructed
by P.Ozsv\'ath and Z.Szab\'o.
In \cite{Gorsky-Namethi} (Theorem 1.4.1 and Corollary 1.5.3) it was shown that (up to a simple change of variables)
the generalized Poincar\'e series of an algebraic link coincides with the
Poincar\'e polynomial (the generating series) of the Heegaard-Floer link homology of it.

\newpage

Addresses:

A. Campillo and F. Delgado:
IMUVA (Instituto de Investigaci\'on en
Matem\'aticas), Universidad de Valladolid.
Paseo de Bel\'en, 7. 47011 Valladolid, Spain.
\newline E-mail: campillo\symbol{'100}agt.uva.es, fdelgado\symbol{'100}agt.uva.es

S.M. Gusein-Zade:
Moscow State University, Faculty of Mathematics and Mechanics, Moscow, GSP-1, 119991, Russia.
\newline E-mail: sabir\symbol{'100}mccme.ru

\end{document}